\documentclass{article}
\usepackage{spconf,amsmath,graphicx}
\usepackage{amsfonts,amssymb}
\usepackage{algorithm}
\usepackage[end]{algpseudocode}
\usepackage{array}

\newcommand{\commnt}[1] {$//$ \textsc{#1} }
\newtheorem{theorem}{\textbf{Theorem}}

\newtheorem{lemma}{\textbf{Lemma}}

\newcommand{\qed}{\begin{flushright} $\Box$ \end{flushright}}

\title{The Non-Bayesian Restless Multi-Armed Bandit: a Case of Near-Logarithmic Regret}
\name{Wenhan Dai$^\dag$\sthanks{Wenhan Dai performed the work
described in this paper during a summer research internship at the
University of Southern California in 2010.}, Yi Gai$^\ddag$, Bhaskar
Krishnamachari$^\ddag$, Qing Zhao$^\S$}
\address{$^\dag$ School of Information Science and Technology, Tsinghua University, China\\ $^\ddag$ Ming Hsieh Department of Electrical Engineering, University of Southern California, U.S.A., \\ $^\S$ Department of Electrical and Computer Engineering, University of California, Davis, U.S.A}
\begin{document}
\ninept
\maketitle
\begin{abstract}
In the classic Bayesian restless multi-armed bandit (RMAB) problem,
there are $N$ arms, with rewards on all arms evolving at each time
as Markov chains with known parameters. A player seeks to activate
$K \geq 1$ arms at each time in order to maximize the expected total
reward obtained over multiple plays. RMAB is a challenging problem
that is known to be PSPACE-hard in general. We consider in this work
the even harder non-Bayesian RMAB, in which the parameters of the
Markov chain are assumed to be unknown \emph{a priori}. We develop
an original approach to this problem that is applicable when the
corresponding Bayesian problem has the structure that, depending on
the known parameter values, the optimal solution is one of a
prescribed finite set of policies. In such settings, we propose to
learn the optimal policy for the non-Bayesian RMAB by employing a
suitable meta-policy which treats each policy from this finite set
as an arm in a different non-Bayesian multi-armed bandit problem for
which a single-arm selection policy is optimal. We demonstrate this
approach by developing a novel sensing policy for opportunistic
spectrum access over unknown dynamic channels. We prove that our
policy achieves near-logarithmic regret (the difference in expected
reward compared to a model-aware genie), which leads to the same
average reward that can be achieved by the optimal policy under a
known model. This is the first such result in the literature for a
non-Bayesian RMAB.
\end{abstract}


\section{Introduction}
\label{sec:intro}

Multi-armed bandit (MAB) problems are fundamental tools for optimal
decision making in dynamic, uncertain environments. In a multi-armed
bandit problem, there are $N$ arms each generating stochastic
rewards, and a player seeks a policy to activate $K \geq 1$ arms at
each time in order to maximize the expected total reward obtained
over multiple plays. MAB problems can be broadly classified as
Bayesian (if player knows the statistical model/parameters of the
reward process for each arm) or non-Bayesian (if the model for the
reward process is a priori unknown to the user). In the case of
non-Bayesian MAB problems, the objective is to design an arm
selection policy that minimizes regret, defined as the gap between
the expected reward that can be achieved by a genie that knows the
parameters, and that obtained by the given policy. It is desirable
to have a regret that grows as slowly as possible over time (if the
regret is sub-linear, the average regret per slot tends to zero over
time, and the policy achieves the maximum average reward that can be achieved under a known model).

A particularly challenging variant of these problems is the restless
multi-armed bandit problem (RMAB)~\cite{Whittle}, in which the
rewards on all arms evolve at each time as Markov chains. Even in
the Baysian case, where the parameters of the Markov chains are
known, this problem is difficult to solve, and has been proved to be
PSPACE hard~\cite{Papadimitriou}. One approach to this problem has
been Whittle's index, which is asymptotically optimal under certain
regimes; however it does not always exist, and even when it does, it
is not easy to compute. It is only in very recent work that
non-trivial tractable classes of RMAB where Whittle's index exists
and is computable have been identified~\cite{LiuZhaoIndex}.

We consider in this work the even harder non-Bayesian RMAB, in which
the parameters of the Markov chain are further assumed to be unknown
\emph{a priori}. Our main contribution in this work is a novel
approach to this problem that is applicable whenever the
corresponding Bayesian RMAB problem has the structure that the
parameter space can be partitioned into a finite number of sets, for
each of which there is a single optimal policy. Our approach
essentially develops a meta-policy that treats these policies as
arms in a different non-Bayesian multi-armed bandit problem for
which a single arm selection policy is optimal for the genie, and
tries to learn which policy from this finite set gives the best
performance.

We demonstrate our approach on a practical problem pertaining to
dynamic spectrum sensing. In this problem, we consider a scenario
where a secondary user must select one of $N$ channels to sense at
each time to maximize its expected reward from transmission
opportunities. If the primary user occupancy on each channel is
modeled to be an identical but independent Markov chain with unknown
parameters, we obtain a non-Bayesian RMAB with the requisite
structure. We develop an efficient new multi-channel cognitive
sensing algorithm for unknown dynamic channels based on the above
approach. We prove for $N=2,3$ that this algorithm achieves regret
(the gap between the expected optimal reward obtained by a
model-aware genie and that obtained by the given policy) that is
bounded uniformly over time $n$ by a function that grows as
$O(G(n)\cdot \log n)$, where $G(n)$ can be any arbitrarily slowly
diverging non-decreasing sequence. This is the first non-Bayesian
RMAB policy that achieves the maximum average reward defined by the
optimal policy under a known model.

There are two parallel investigations on non-Bayesian RMAB problems
given in \cite{TekinLiu,LiuLiuZhao}, where a more general RMAB model
is considered but under a much weaker definition of regret.
Specifically, in \cite{TekinLiu,LiuLiuZhao}, regret is defined with
respect to the maximum reward that can be offered by a \emph{single
arm/channel}. Note that for RMAB with a known model, staying with
the best arm is suboptimal. Thus, a sublinear regret under this
definition does not imply the maximum average reward, and the
deviation from the maximum average reward can be arbitrarily large.

\section{A New Approach for non-Bayesian RMAB}

We first describe a structured class of finite-option Bayesian RMAB
problems that we will refer to as $\Psi_m$. Let $\mathcal{B}(P)$ be
a Bayesian RMAB problem with the Markovian evolution of arms
described by the transition matrix $P$. We say that $\mathcal{B}(P)
\in \Psi_m$ if and only if there exists a partition of the parameter
values $P$ into a finite number of $m$ sets $\{ S_1, S_2, ... S_m\}
$ and a set of policies $\pi_i$ $(\forall i = 1\ldots m)$ that do
not assume knowledge of $P$ and are optimal whenever $P \in S_i$.
Despite the general hardness of the RMAB problem, problems with such
structure do indeed exist, as has been shown
in~\cite{zhao:krishnamachari:twc2008, ahmad:liu:it2009,
LiuZhaoIndex}.

We propose a solution to the non-Bayesian version of the problem
that leverages the finite solution option structure when we have
that the corresponding Bayesian version $\mathcal{B}(P) \in \Psi_m$.
In this case, although the player does not know the exact parameter
$P$, it must be true that one of the $m$ policies $\pi_i$ will yield
the highest expected reward (corresponding to the set $S_i$ that
contains the true, unknown $P$). These policies can thus be treated
as arms in a different non-Bayesian multi-armed bandit problem for
which a single-arm selection policy is optimal for the genie. Then,
a suitable meta-policy that sequentially operates these policies
while trying to minimize regret can be adopted. This can be done
with an algorithm based on the well-known schemes proposed by Lai
and Robbins~\cite{Lai:Robbins}, and Auer \emph{et
al}~\cite{Auer:2002}.

One subtle issue that must be handled in adopting such an algorithm
as a meta-policy is how long to play each policy. An ideal constant
length of play could be determined only with knowledge of the
underlying unknown parameters $P$, so our approach is to have the
duration for which each policy is operated slowly increase over
time.

In the following, we demonstrate this novel meta-policy approach to
the dynamic spectrum access problem discussed
in~\cite{zhao:krishnamachari:twc2008, ahmad:liu:it2009} where the
Bayesian version of the RMAB has been shown to belong to the class
$\Psi_2$. For this problem, we show that our approach yields an
algorithm with provably near-logarithmic regret, thus achieving the
same average reward offered by the optimal RMAB policy under a known
model.

\section{A Dynamic Spectrum Access Problem}
\label{sec:problem}

We consider a slotted system where a secondary user is trying to
access $N$ independent channels, with the availability of each
channel evolving as a two-state Markov chain with identical
transition matrix $\mathbf{P}$ that is \emph{a priori} unknown to
the user. The user can only see the state of the sensed channel. If
the user selects channel $i$ at time $t$, and upon sensing finds the
state of the channel $S_i(t)$ to be 1, it receives a unit reward for
transmitting. If it instead finds the channel to be busy, i.e.,
$S_i(t) = 0$, it gets no reward at that time. The user aims to
maximize its expected total reward (throughput) over some time
horizon by choosing judiciously a sensing policy that governs the
channel selection in each slot. We are interested in designing
policies that perform well with respect to \emph{regret}, which is
defined as the difference between the expected reward that could be
obtained using the omniscient policy $\pi^*$ that knows the
transition matrix $\mathbf{P}$, and that obtained by the given
policy $\pi$. The regret at time $n$ can be expressed as:
\begin{equation}
\begin{split}
 R(\mathbf{P},\Omega(1),n) & =
 E^{\pi^*}[\Sigma_{t=1}^{n}Y^{\pi^*}(\mathbf{P},\Omega(1),t)]
 \\
 & -
 E^{\pi}[\Sigma_{t=1}^{n}Y^{\pi}(\mathbf{P},\Omega(1),t)],
\end{split}
\end{equation}
where $\omega_i$ is the initial probability that $S_i(1)=1$,
$\mathbf{P}$ is the transition matrix of each channel,
$Y_{\pi^*}(\mathbf{P},\Omega(1),t)] $ is the reward obtained in time
$t$ with the optimal policy, $ Y_{\pi}(\mathbf{P},\Omega(1),t)$ is
the reward obtained in time $t$ with the given policy. We denote
$\Omega(t)\triangleq [\omega_1(t), \ldots, \omega_N(t)]$ as the
belief vector where $\omega_i(t)$ is the conditional probability
that $S_i(t) = 1$ (and let
$\Omega(1)=[\omega_1(1),\ldots,\omega_N(1)]$ denote the initial
belief vector used in myopic sensing
algorithm~\cite{zhao:krishnamachari:twc2008}).

\section{Sensing Unknown Dynamic Channels}
\label{sec:algorithm}

As has been shown in~\cite{zhao:krishnamachari:twc2008}, the myopic policy has a
simple structure for switching between channels that depends only on
the correlation sign of the transition matrix $\mathbf{P}$, i.e.
whether $p_{11}\geq p_{01}$ (positively correlated) or $p_{11}<
p_{01}$ (negatively correlated).

In particular, if the channel is positively correlated, then the
myopic policy corresponds to
\begin{itemize}
\item \textbf{Policy $\pi_1$}: stay on a channel whenever it shows a
``1''
and switch on a ``0'' to the channel visited the longest ago.
\end{itemize}

If the channel is negatively correlated, then it corresponds to
\begin{itemize}
\item \textbf{Policy $\pi_2$:} staying on a channel when it shows a ``0'', and
switching as soon as ``1'' is observed, to either the channel most
recently visited among those visited an even number of steps before,
or if there are no such channels, to the one visited the longest
ago.
\end{itemize}

Furthermore, it has been shown in~\cite{zhao:krishnamachari:twc2008,
ahmad:liu:it2009} that the myopic policy is optimal for $N=2,3$, and
for any $N$ in the case of positively correlated channels (the
optimality of the myopic policy for $N>3$ negatively correlated
channels is conjectured for the infinite-horizon case). As a
consequence, this special class of RMAB has the required finite
dependence on its model as described in Sec.~2; specifically, it
belongs to $\Psi_2$. We can thus apply the general approach based on
the concept of meta-policy. Specifically, the algorithm treats these
two policies as arms in a classic non-Bayesian multi-armed bandit
problem, with the goal of learning which one gives the higher
reward.

A key question is how long to operate each arm at each step. It
turns out from the analysis we present in the next section that it
is desirable to slowly increase the duration of each step using any
(arbitrarily slowly) divergent non-decreasing sequence of positive
integers $\{K_n\}_{n=1}^\infty$.

The channel sensing policy we thus construct is shown in Algorithm
\ref{alg:sensing}.

\begin{algorithm} [ht]
\caption{Sensing Policy for Unknown Dynamic Channels}
\label{alg:sensing}
\begin{algorithmic}[1]
\State \commnt{ Initialization}

\State Play policy $\pi_1$ for $K_1$ times, denote $\hat{A}_1$ as
the sample mean of these $K_1$ rewards \State Play policy $\pi_2$
for $K_2$ times, denote $\hat{A}_2$ as the sample mean of these
$K_2$ rewards \State $\hat{X}_1 = \hat{A}_1$, $\hat{X}_2 =
\hat{A}_2$ \State $n = K_1+K_2$ \State $i=3$, $i_1=1$, $i_2=1$

\State \commnt{Main loop}

\While {1}
    \State Find $j$ such that
    $j = \arg\max\frac{\hat{X}_j}{i_j}+\sqrt{\frac{3\ln{n}}{i_j}}$ \label{line:9}
    \State $i_j = i_j+1$
    \State Play policy $\pi_j$ for $K_i$ times, let $\hat{A}_{j}(i_j)$ record the sample mean of these $K_i$ rewards
    \State $\hat{X}_j = \hat{X}_j + \hat{A}_{j}(i_j)$
    \State $i = i+1$
    \State $n = n + K_i$;
\EndWhile
\end{algorithmic}
\end{algorithm}

\section{Regret Analysis}
\label{sec:regret}

We first define the discrete function $G(n)$ which represents the
value of $K_i$ at the $n^{th}$ time step in
Algorithm~\ref{alg:sensing}:
\begin{equation}G(n) = \min\limits_I K_I~~s.t.~~\sum\limits_{i=1}^I K_i \geq  n
\end{equation}
Note that since $K_i$ can be any arbitrarily slow non-decreasing
diverging sequence $G(n)$ can also grow arbitrarily slowly.

The following theorem states that the regret of our algorithm grows
close to logarithmically with time.

\begin{theorem} \label{theorem:regret} For the dynamic spectrum access
problem with $N=2,3$ i.i.d. channels with unknown transition matrix
$\mathbf{P}$, the expected regret with Algorithm \ref{alg:sensing}
after $n$ time steps is at most $Z_1 G(n)\ln(n) + Z_2 \ln(n) + Z_3
G(n) + Z_4$, where $Z_1,Z_2,Z_3,Z_4$ are constants only related to
$\mathbf{P}$.
\end{theorem}

The proof of Theorem \ref{theorem:regret}, presented in the
appendix, uses two interesting lemmas we have developed that we
present here without proof. The first lemma is a non-trivial variant
of the Chernoff-Hoeffding bound, that allows for bounded differences
between the conditional expectations of sequence of random variables
that are revealed sequentially:


\begin{lemma} \label{lemma:chernoff} Let $X_1,\cdots,X_n$ be random
variables with range $[0,b]$ and such that
$|E[X_t|X_1,\cdots,X_{t-1}]- \mu | \leq C$. $C$ is a constant number
such that $0< C < \mu$. Let $S_n = X_1+\cdots+X_n$. Then for all  $a
\geq 0$,
\begin{equation} \label{eqn:chernoff1}
P\{S_n \geq n(\mu+C) + a\}\leq e^{-2(\frac{a(\mu-C)}{b(\mu+C)})^2/n}
\end{equation}
\begin{equation} \label{eqn:chernoff2}
P\{S_n \leq n(\mu-C) - a\}\leq e^{-2(a/b)^2/n}
\end{equation}
\end{lemma}




The second lemma states that the expected loss of reward for either
policy due to starting with an arbitrary initial belief vector
compared to the reward $U_i(P)$ that would obtained by playing the
policy at steady state is bounded by a constant $C_i(P)$ that
depends only on the policy used and the transition matrix. These
constants can be calculated explicitly, but we omit the details for
brevity.

\begin{lemma} \label{lemma:2} For any initial belief vector $\Omega(1)$
and any positive integer $L$, if we use \emph{policy $\pi_i$}
($i=1,2$) for $L$ times, and the summed expectation of the rewards
for these $L$ steps is denoted as
$E^{\pi_i}[\Sigma_{t=1}^{L}Y^{\pi_i}(\mathbf{P},\Omega(1),t)]$, then
\begin{equation}
|E^{\pi_i}[\Sigma_{t=1}^{L}Y^{\pi_i}(\mathbf{P},\Omega(1),t)]-L
\cdot U_i(P)|<C_i(P)
\end{equation}
\end{lemma}

\textbf{Remark:} Theorem 1 has been stated for the cases $N=2,3$,
which are the only cases when the Myopic policy has been proved to
be optimal for the known-parameter case for all values of $P$. In
fact, our proof shows something even stronger than this: that
Algorithm \ref{alg:sensing} yields the claimed near-logarithmic
regret with respect to the Myopic policy for any $N$. The Myopic
policy is known to be always optimal for $N=2,3$, and for any $N$ so
long as the Markov chain is positively correlated. In case it is
negatively correlated, it is an open question whether it is optimal
for an infinite horizon case. If this were to be true, the algorithm
we have presented would also offer near-logarithmic regret
asymptotically as the time variable $n$ increases, for any $N$.




\begin{thebibliography}{1}

\bibitem{Whittle} P. Whittle, ``Restless Bandits: Activity Allocation in a Changing World," Journal of Applied Probability, Vol. 25,
1988.

\bibitem{Papadimitriou} C. H. Papadimitriou and J. N. Tsitsiklis, ``The Complexity Of Optimal Queueing Network Control," \emph{Mathematics of Operations Research}, Vol. 24, 1994.

\bibitem{LiuZhaoIndex}
K.~Liu and Q.~Zhao, ``"Indexability of restless
bandit problems and optimality of Whittle index for dynamic multichannel
access,'' \emph{IEEE Trans. Inf. Theory}, vol. 56, no. 11, November, 2010.

\bibitem{zhao:krishnamachari:twc2008}
Q.~Zhao, B.~Krishnamachari, and K.~Liu, ``On myopic sensing for
multi-channel opportunistic access: structure, optimality, and
performance,'' \emph{IEEE Transactions on Wireless Communications},
2008.

\bibitem{ahmad:liu:it2009}
S.~Ahmad, M.~Liu, T.~Javidi, Q.~Zhao, and B.~Krishnamachari,
``Optimality of myopic sensing in multi-channel opportunistic
access,'' \emph{IEEE Transactions on Information Theory}, 2009.

\bibitem{Lai:Robbins}
T.~Lai and H.~Robbins, ``Asymptotically efficient adaptive
allocation rules,'' \emph{Advances in Applied Mathematics}, vol. 6,
no. 1, 1985.

\bibitem{Auer:2002}
P.~Auer, N.~Cesa-Bianchi, and P.~Fischer, ``Finite-time analysis of
the multiarmed bandit problem,'' \emph{Machine Learning}, 47(2-3),
2002.


\bibitem{TekinLiu} C. Tekin and M. Liu, ``Online Learning in Opportunistic Spectrum Access: A Restless
Bandit Approach," Arxiv pre-print http://arxiv.org/abs/1010.0056,
October 2010.

\bibitem{LiuLiuZhao}
 H. Liu, K. Liu and Q. Zhao, ``Logrithmic Weak Regret of Non-Bayesian Restless Multi-Armed Bandit,''
submitted to \emph{ICASSP}, October, 2010.

\end{thebibliography}

\section{Appendix}

%
%


\textbf{Proof of Theorem \ref{theorem:regret}}

We first derive a bound on the regret for the case when
$p_{01}<p_{11}$. In this case, policy $\pi_1$ would be the optimal.
Based on Lemma \ref{lemma:2}, the difference of
$E^{\pi_1}[\Sigma_{t=1}^{n}Y^{\pi_1}(\mathbf{P},\Omega(1),t)]$ and
$U_1\cdot n$ is no more than $C_1$, therefore we only need to prove
$R'(\mathbf{P},\Omega(1),n)$, the regret in the case when policy
$\pi_1$ is optimal, which is defined as $R'(\mathbf{P},\Omega(1),n)
\triangleq U_1 \cdot
n-E^{\pi_1}[\Sigma_{t=1}^{n}Y_{\pi_1}(\mathbf{P},\Omega(1),t)]$, is
at most $Z_1 G(n)\ln(n) + Z_2 \ln(n) + Z_3 G(n) + Z_4$,
$Z_1,Z_2,Z_3,Z_4$ are constants only related to $\mathbf{P}$.

The regret comes from two parts: the regret when using policy
$\pi_2$; the regret between $U_1$ and
$E^{\pi_1}[Y_{\pi_1}(\mathbf{P},\Omega(1),t)]$ when using policy
$\pi_1$. From Lemma \ref{lemma:2}, we know that each time when we
switch from policy $\pi_2$ to policy $\pi_1$, at most we lose a
constant-level value from the second part. So if the number of
selections of policy $\pi_2$ in Line \ref{line:9} of Algorithm
\ref{alg:sensing} is bounded by $O(\ln{n})$, both parts of the
regret can be bounded by $O(G(n)\cdot\ln{n})$.

For ease of exposition, we discuss the slots $n$ such that $G||n$ ,
where $G||n$ denotes that time $n$ is the end of successive $G(n)$
plays.

We define $q$ as the smallest index such that $K_q \geq \lceil
\frac{C_1+C_2}{|U_1-U_2|} \rceil$. Note that it is possible to
define $\alpha(U_1, C_1, \mathbf{P})$ such that if policy $\pi_1$ is
played $s_1 > \alpha$ times,
\begin{equation} \label{eqn:p1}
\exp(-2(q(U_1-\frac{C_1}{K_q})-s\sqrt{\frac{3\ln{t}}{s}})^2/(s-q))
\leq 2t^{-4}
\end{equation}
We could also define $\beta(U_2, C_2, \mathbf{P})$ such that if
policy $\pi_2$ is played $s_2 > \beta$ times,
\begin{equation}\label{eqn:p2}
\exp(\frac{-2(q\frac{U_2-C_2/K_q}{U_2+C_2/K_q}(U_2+\frac{C_2}{K_q}-1)+s\sqrt{\frac{3\ln{t}}{s}})^2}{s-q})\\
\leq 2t^{-4}
\end{equation}
Moreover, there exists $\gamma = \lceil
\max\{{5\alpha+1,e^{4\alpha/3},5\beta+1,e^{4\beta/3}}\}\rceil$ such
that when $G(n)>K_{\gamma}$, policy $\pi_1$ is played at least
$\alpha$ times and policy $\pi_2$ is played at least $\beta$ times.

Denote $T(n)$ as the number of times we select policy $\pi_2$ up to
time $n$. Then, for any positive integer $l$, we have:
\begin{equation}
\begin{split} T(n) & =
1+\sum_{t=K_1+K_2,G||t}^{n}\mathbb{I}\{\frac{\hat{X}_1(t)}{i_1(t)}+\sqrt{\frac{3\ln{t}}{i_1(t)}}
\\
 & \quad <
\frac{\hat{X}_2(t)}{i_2(t)}+\sqrt{\frac{3\ln{t}}{i_2(t)}}\} \\
& \leq l+\gamma+\\
&\sum_{t=K_1+\cdots+K_\gamma,G||t}^{n}\sum_{s_1=\alpha}^{\alpha(t),t=K_1+\cdots+K_{\alpha(t)}}\sum_{s_2=max(\beta,l)}^{\beta(t),t=K_1+\cdots+K_{\beta(t)}}\\
&\mathbb{I}\{\frac{\hat{X}_{1,s_1}}{s_1}+\sqrt{\frac{3\ln{t}}{s_1}}\leq\frac{\hat{X}_{2,s_2}}{s_2}+\sqrt{\frac{3\ln{t}}{s_2}}\}
\end{split}
\end{equation}
where $\mathbb{I}\{x\}$ is the index function defined to be 1 when
the predicate $x$ is true, and 0 when it is false predicate;
$i_j(t)$ is the number of times we select policy $\pi_j$ when up to
time $t, \forall j=1,2$; $\hat{X}_j(t)$ is the sum of every sample
mean for $K_i$ plays up to time $t$; $\hat{X}_{i,s_i}$ is the sum of
every sample mean for $K_{s_i}$ times using policy $\pi_i$.

The condition
$\{\frac{\hat{X}_{1,s_1}}{s_1}+\sqrt{\frac{3\ln{t}}{s_1}}\leq\frac{\hat{X}_{2,s_2}}{s_2}+\sqrt{\frac{3\ln{t}}{s_2}}\}$
implies that at least one of the following must hold:
\begin{equation} \label{eqn:inequ1}
\frac{\hat{X}_{1,s_1}}{s_1}\leq
U_1-\frac{C_1}{K_q}-\sqrt{\frac{3\ln{t}}{s_1}}
\end{equation}
\begin{equation} \label{eqn:inequ2}
\frac{\hat{X}_{2,s_2}}{s_2}\geq
U_2+\frac{C_2}{K_q}+\frac{U_2+C_2/K_q}{U_2-C_2/K_q}\sqrt{\frac{3\ln{t}}{s_2}}
\end{equation}
\begin{equation} \label{eqn:inequ3}
U_1-\frac{C_1}{K_q}<U_2+\frac{C_2}{K_q}+(1+\frac{U_2+C_2/K_q}{U_2-C_2/K_q})\sqrt{\frac{3\ln{t}}{s_2}}
\end{equation}

Note that
$\hat{X}_{1,s_1}=\hat{A}_{1,1}+\hat{A}_{1,2}+\cdots+\hat{A}_{1,s_1}$,
where $\hat{A}_{1,i}$ is sample average reward for the $i_{th}$ time
selecting policy $\pi_1$. Due to the definition of $\alpha$ and
$K_q$, the expected value of $\hat{A}_{1,i}$ is between
$U_1-\frac{C_1}{K_q}$ and $U_1+\frac{C_1}{K_q}$ if $i \geq q$ (Lemma
\ref{lemma:2}). Then applying Lemma \ref{lemma:chernoff}, and the
results in (\ref{eqn:p1}) and (\ref{eqn:p2}),
\begin{equation}
\begin{split}
&Pr(\frac{\hat{X}_{1,s_1}}{s_1}\leq
U_1-\frac{C_1}{K_q}-\sqrt{\frac{3\ln{t}}{s_1}}) \leq 2t^{-4}
\end{split}
\end{equation}
\begin{equation} Pr(\frac{\hat{X}_{2,s_2}}{s_2}\geq
U_2+\frac{C_2}{K_q}+\frac{U_2+C_2/K_q}{U_2-C_2/K_q}\sqrt{\frac{3\ln{t}}{s_2}})
\leq 2t^{-4}\end{equation}


For
$\lambda(n)=\lceil(3(1+\frac{U_2+C_2/K_q}{U_2-C_2/K_q})^2\ln{n})/(U_1-U_2-\frac{C_1+C_2}{K_q})^2\rceil$,
(\ref{eqn:inequ3}) is false. So we get:
\begin{equation}
\begin{split}
&E(T(n)) \leq
\lambda(n) +\gamma+
\Sigma_{t=1}^{\infty}\Sigma_{s_1=1}^{t}\Sigma_{s_2=1}^{t}4t^{-4}\\
&\leq
\lambda(n)+\gamma+\frac{2\pi^2}{3}.
\end{split}
\end{equation}

Therefore, we have:
\begin{equation}
\begin{split}
& R'(\mathbf{P},\Omega(1),n) \leq G(n)+ \\
& \quad\quad  ((U_1-U_2)G(n)+C_1)
(\lambda(n)+\gamma+\frac{2\pi^2}{3})
\end{split}
\end{equation}

This concludes the bound in case $p_{11}>p_{01}$. The derivation of
the bound is similar for the case when $p_{11}\leq p_{01}$ with the
key difference of $ \gamma'$ instead of $\gamma$, and the $C_1, U_1$
terms being replaced by $C_2, U_2$ and vice versa. Then we have that
the regret in either case has the following bound:

\begin{equation}
\begin{split}
 R(\mathbf{P},\Omega(1),n) & \leq G(n)+
(|U_2-U_1|G(n)+\max\{C_1,C_2\})\\
&
(\frac{3(1+\max\{\frac{U_1+C_1/K_q}{U_1-C_1/K_q},\frac{U_2+C_2/K_q}{U_2-C_2/K_q}\})^2\ln{n}}{(|U_2-U_1|-\frac{C_1+C_2}{K_q})^2}
\\
& + 1 +{\max\{\gamma,\gamma'}\}+\frac{2\pi^2}{3})
\end{split}
\end{equation}

This inequality can be readily translated to the simplified form of
the bound given in the statement of Theorem 1, where:
\begin{equation}
\begin{split}
& Z_1 = |U_2-U_1|
\frac{3(1+\max\{\frac{U_1+C_1/K_q}{U_1-C_1/K_q},\frac{U_2+C_2/K_q}{U_2-C_2/K_q}\})^2\ln{n}}{(|U_2-U_1|-\frac{C_1+C_2}{K_q})^2}\\
& Z_2 = \max\{C_1,C_2\}
\frac{3(1+\max\{\frac{U_1+C_1/K_q}{U_1-C_1/K_q},\frac{U_2+C_2/K_q}{U_2-C_2/K_q}\})^2\ln{n}}{(|U_2-U_1|-\frac{C_1+C_2}{K_q})^2}\\
& Z_3 = |U_2-U_1| (1 +{\max\{\gamma,\gamma'}\}+\frac{2\pi^2}{3}) +
1\\ \nonumber
& Z_4 = \max\{C_1,C_2\} (1
+{\max\{\gamma,\gamma'}\}+\frac{2\pi^2}{3})
\end{split}
\end{equation}
\qed

\end{document}